% Plain TeX file of the paper Meromorphic extendibility and the argument principle
% Written by Josip Globevnik
% Final version prepared on December 1, 2006
\magnification 1200
\def\R{{\rm I\kern-0.2em R\kern0.2em \kern-0.2em}}
\def\N{{\rm I\kern-0.2em N\kern0.2em \kern-0.2em}}
\def\P{{\rm I\kern-0.2em P\kern0.2em \kern-0.2em}}
\def\B{{\rm I\kern-0.2em B\kern0.2em \kern-0.2em}}
\def\Z{{\rm I\kern-0.2em Z\kern0.2em \kern-0.2em}}
\def\C{{\bf \rm C}\kern-.4em {\vrule height1.4ex width.08em depth-.04ex}\;}

\def\D{{\Delta}}

\def\cW{{\cal W}}

\font\ninerm=cmr8
\ 
\vskip 22mm
\centerline {\bf MEROMORPHIC EXTENDIBILITY AND THE }

\centerline {\bf ARGUMENT PRINCIPLE}
\vskip 8mm
\centerline{Josip Globevnik}
\vskip 8mm
{\noindent \ninerm ABSTRACT \ \ Let $\D $ be the open unit disc in $\C$. Given a 
continuous function $\varphi \colon \ 
b\D \rightarrow \C\setminus \{ 0\}$
 denote by $\cW (\varphi )$ the winding number of $\varphi$ around the origin. We prove that 
a continuous function $f\colon\ b\D\rightarrow \C$ extends meromorphically through 
$\D $ if and only if there is a 
number $N\in \N\cup\{ 0\}$ such 
that $\cW (Pf+Q)\geq -N$  for every pair $P, Q$ of polynomials such 
that $Pf+Q\not= 0$ on $b\D$. If this is the case then 
the meromorphic extension has at most $N$ poles in $\D$.}
\vskip 6mm
\bf 1.\ Introduction and the main result \rm 
\vskip 2mm
Let $\D $ be the open unit disc in $\C$ and let $f\colon \ b\D\rightarrow \C$ be a continuous function.
We say 
that $f$  extends holomorphically  through $\D$ if $f$ admits a continous extension $\tilde f$ 
to $\overline \D$ which is holomorphic on $\D$. If this 
is the case 
then we say that $f$ (or $\tilde f$) belongs to the disc algebra. Denote by 
$\overline \C =\C\cup\{ \infty\}$ the Riemann sphere.
We say that $f$ extends meromorphically through $\D$ if 
there is a finite set $A\subset \D$ such that $f$ has a continuous 
extension to $\overline \D\setminus A$ which is holomorphic on $\D\setminus A$
and has a pole at each point of $A$. Equivalently, $f$ extends meromorphically 
through $\D $ if it has a continous extension 
$\tilde f\colon \overline\D \rightarrow \overline \C$ which, 
as a function to $\overline \C$,  is holomorphic on $\D$. 

Given a continuous function  $\varphi\colon \ b\D\rightarrow \C\setminus \{ 0\}$ we denote 
by $\cW (\varphi )$ the winding number of
$\varphi $ (around the origin). So $\cW (\varphi )$ equals $1/(2\pi )$times the change of argument of $
\varphi (z )$ as $z$ runs once around $b\D $ counterclockwise. 

In the present paper we show that meromorphic extendibility can be characterized 
in terms of the argument principle. For holomorphic extendibility this is already known:
\vskip 2mm
\noindent \bf THEOREM 1.0\ \rm[G2] \it A continuous function $f\colon\ b\D\rightarrow \C$ 
extends holomorphically through $b\D$ if and only if $\cW (f+Q)\geq 0 $ for every polynomial $Q$ such that
$f+Q\not= 0$ on $b\D$. \rm
\vskip 2mm

If a continuous function $f\colon \ b\D \rightarrow \C\setminus \{ 0\} $ 
extends meromorphically through $\D $ then $\cW (f)\geq -N$ where $N$ is the number
of poles of the meromorphic extension $\tilde f$ (counted with multiplicity). 
Indeed, by the argument principle, 
$$
\cW (f) = \nu _0(\tilde f)-\nu _p(\tilde f) = \nu _0 (\tilde f) - N\geq -N
$$
where $\nu _0 (\tilde f)$ is the number of zeros of $\tilde f$ on $\D $ 
and $\nu _p(\tilde f)$ is the number of poles of 
$\tilde f$ on $\D$.

Let $f \colon \ b\D\rightarrow \C$ be a continuous function on $b\D$ which extends meromorphically
through $\D $ and whose meromorphic extension $\tilde f$ has 
$N$ poles on $\D $. Then
$
\cW (Pf+Q) \geq -N 
$
for all polynomials $P, Q$ such that $Pf+Q\not = 0$ on $b\D$. Indeed, $P\tilde f + Q$, the 
meromorphic extension of $Pf+Q$, has no other poles than $\tilde f$ and therefore, 
by the argument principle, \ $\cW (Pf+Q)\geq -N$. 
The following theorem, our main result, tells 
that this property characterizes meromorphic extendibility. 
\vskip 2mm
\noindent \bf THEOREM 1.1\ \ \it A continuous function $f\colon \ b\D\rightarrow \C$ 
extends meromorphically through 
$\D $ if and only if there is an $N\in\N\cup \{0 \} $ such that 
$$
\cW (Pf+Q)\geq -N
\eqno (1.1)
$$
for all polynomials $P, Q$ such that $Pf+Q\not= 0 $ on $b\D $. If this 
is the case then the meromorphic extension of $f$ has at most $N$ poles 
in $\D$, counting multiplicity. 
\vskip 4mm
\bf 2. Fourier series \rm
\vskip 2mm
In this section we recall some well known facts. 
\vskip 2mm
Let f be a continuous function on $b\D $. For each integer $n$ let 
$$\hat f (n) = {1\over{2\pi}}\int _{-\pi}^\pi e^{-i n \theta}f(e^{i\theta }) d\theta 
$$
so that
$$
\sum_{n=-\infty}^\infty \hat f(n)e^{in\theta }
$$
is the Fourier series of $f$. We have
$$
\sum_{n=-\infty}^\infty |\hat f(n)|^2 <\infty .
\eqno (2.1)
$$
Define the functions $F$ and $G$ by
$$
F(z)= \hat f (0) +\hat f (1) z + \hat f(2) z^2 +\cdots \ \ \ \ \ (z\in\D)
\eqno 
$$
$$
G(z) = \overline{\hat f(-1)}z +\overline 
{\hat f(-2)}z^2+\cdots \ \ \ (z\in\D).
$$
The functions $F$ and $G$ are holomorphic on $\D$ and by (2.1) they belong 
to the space $H^2$ [R]. 

The function $f$ belongs to the disc algebra if and only if $\hat f(n)=0$ for all $n<0$ or, 
equivalently, if and only if $G\equiv 0$.

Suppose now that $f$ is smooth. Then the Fourier series converges uniformly 
to $f$. The functions $F$ and $G$ belong to the disc algebra and have smooth boundary values. 
We have
 $$  
f(z) = F(z) + \overline {G(z)}\ \ (z\in b\D).
$$
We shall need the following
\vskip 2mm
\noindent \bf PROPOSITION 2.1\ \it Let $\Phi \colon\ b\D\rightarrow \C$ be a continuous function.
Given $N\in \N$ there is a nonzero polynomial $P$ of degree not exceeding $N$ such that 
$\hat {(P\Phi )}(j)=0\ (-N\leq j\leq -1)$ . \rm
\vskip 2mm
\noindent\bf Proof.\ \rm $P\Phi $ is continuous on $\D $ and a direct computation shows that for 
each integer $j$ we have 
$$
\hat{(P\Phi)}(j) = \hat P(0)\hat\Phi (j) + \hat P(1)\hat \Phi (j-1) +\cdots + \hat P (N)\hat\Phi (j-N)
$$
so $\hat{(PF)}(j) = 0\ (-N\leq j\leq -1)$ gives the homogeneous system 
$$
\eqalign{
&\hat P (0)\hat\Phi (-N) + \hat P(1)\hat\Phi (-N-1)+\cdots + \hat P(N)\hat \Phi (-2N)=0 \cr
&\hat P (0)\hat\Phi (-N+1) + \hat P(1)\hat\Phi (-N)+\cdots + \hat P(N)\hat \Phi (-2N+1)=0 \cr
&\cdots \cr
&\hat P (0)\hat\Phi (-1) + \hat P(1)\hat\Phi (-2)+\cdots + \hat P(N)\hat \Phi (-N-1)=0 \cr}
$$
of $N$ linear equations with $N+1$ unknowns $\hat P(0),\hat P(1),\cdots ,\hat P (N)$ which 
always has a nontrivial solution. This completes the proof.
\vskip 4mm
\bf 3.\ The smooth case \rm
\vskip 2mm
\noindent \bf THEOREM 3.1\ \ \it Suppose that $f\colon\ b\D \rightarrow \C$ is of the form
$$
f(z)= G(z) + \overline {H(z)}\ \ \ (z\in b\D )
$$
where the functions $G$ and $H$ belong to the disc algebra and $H$ has smooth 
boundary values. Assume that $N\in\N\cup\{ 0\}$ and that 
$$
\left.\eqalign{
\cW (Pf+Q)\geq -N\ \hbox{\ whenever\ } P, Q, \hbox{\ are polynomials,}\cr
deg (P)\leq N,\ \hbox{such that \ } Pf+Q\not= 0\hbox{\ on\ } b\D .}\right\} \eqno (3.1)
$$ 
Then $f$ extends meromorphically through $\D$ and the meromorphic 
extension has at most $N$ poles in $\D$, counting multiplicity.
\vskip 2mm
\noindent\bf REMARK 3.2\ \rm To prove Theorem 1.1  we shall later use Theorem 
3.1 only in the special case when 
$H$ is a rational function holomorphic in a neighbourhood of  $\overline \D$ so with this 
in mind, we may assume as much
smoothness as we want. In the proof of Theorem 3.1 below it is
enough to assume 
that $H|b\D$ belongs to the Lipschitz class $C^\alpha $ 
with $\alpha > 1/2$. 

Before proceeding observe that if $f$ is as in Theorem 3.1 and $P$ is a 
polynomial then $Pf$ has the same form. 
Indeed, we have $Pf=PG + P\overline H $ on $b\D $ where the function 
$z\mapsto P(z)\overline{H(z)}$ is smooth on $b\D$ so on $b\D$ we have 
$P\overline H=F_2+\overline {G_1} $ where $F_2, G_1$ belong to the disc 
algebra and have smooth boundary values. So on $b\D $ we have 
$Pf=PG+F_2+\overline {G_1}= F_1 + \overline{G_1} $ 
where $F_1, G_1 $ are in the disc algebra and $G_1$ has smooth boundary values. 
\vskip 2mm
\noindent\bf Proof of Theorem 3.1. \rm Assume that $f$ is as in Theorem 3.1 and that (3.1) holds for 
some $N\in\N\cup\{ 0\}$. If $N=0$ then it is known that $f$ extends 
holomorphically through $\D $ [G2]. Assume that $N\geq 1$. 
By Proposition 2.1 there is a polynomial $P,\ deg(P)\leq N$, such that
$$
\hat{(Pf)}(-1)=\hat{(Pf)}(-2)=\cdots = \hat{(Pf)}(-N) = 0
\eqno (3.2)
$$
Now, $Pf=F_1 +\overline{G_1}$ on $b\D$ where $F_1, G_1 $ are in 
the disc algebra and $G_1$ has smooth boundary values. 
With no loss of generality assume that $G_1(0)=0$. By (3.2) $G_1= z^{N+1}G_2$ where $G_2$ is 
again in the disc algebra and has smooth boundary values so that 
$$
P(z)f(z)-F_1(z) = \overline{z^{N+1}G_2(z)}\ \ (z\in b\D).
$$
Suppose for a moment that $G_2\not\equiv 0$. 
We show that there is a constant $\alpha \in\C$ such that 
$\overline{z^{N+1}G_2(z)}+\alpha \not= 0 \ (z\in b\D )$ and 
$$
\cW (\overline{z^{N+1}G_2(z)} +\alpha )\leq -N-1,
\eqno (3.3)
$$ 
that is,
$$
\cW ({z^{N+1}G_2(z)} +\overline \alpha )\geq N+1.
$$ 
The function $\Phi(z)=z^{N+1}G_2(z)$ belongs to the disc algebra and has smooth boundary values. 
It has zero of order 
at least $N+1$ at the origin. If $\Phi(z)\not=0\ (z\in b\D)$ then put $\alpha = 0$.  
In this case $\cW (\Phi )$ equals the number of zeros of $\Phi$ in $\D $ so $\cW (\Phi ) \geq N+1$. Suppose
now that $\Phi(b\D )$ contains the origin. Since 
$\Phi$ has smooth boundary values it follows that $\Phi(b\D ) $ is nowhere dense.
So there are $\alpha $, arbitrarily close to the origin such that 
$\Phi(z)+\overline\alpha \not=0\ (z\in b\D )$. Let 
$\nu \geq N+1$ be the multiplicity of the zero of $\Phi$ at the origin. A standard use 
of the argument principle on a sufficiently small disc $D$ centered at the origin shows that 
for any $\alpha $ sufficiently 
close to the origin, $\alpha \not= 0$, the function $z\mapsto \Phi(z)+\overline\alpha$ 
has exactly $\nu $ zeros on $D$ which are arbitrarily close to the origin provided 
that $\alpha $ is sufficiently close to the origin. Thus, if $\alpha \not=0$ is sufficiently 
close to $0$ and $\Phi(z)+\overline\alpha \not=0\ (z\in b\D )$ then 
$\Phi+\overline\alpha $ has $\nu $ zeros in a neighbourhood of the origin so the argument 
principle, now applied to the function $\Phi-\overline\alpha$ on $\D $, implies that $\cW (\Phi+\overline \alpha)\geq \nu\geq N+1$ so that 
(3.3) holds. It follows that $\cW (Pf-F_1+\alpha )\leq -N-1$. A sufficiently good polynomial approximation
$Q$ of $-F_1+\alpha $ then satisfies $\cW (Pf+Q)\leq -N-1$, contradicting (3.1). It follows 
that $G_2\equiv 0$ so $Pf=F_1$ on 
$b\D$, that is, $Pf$ belongs to the disc algebra. We need
\vskip 2mm
\noindent\bf PROPOSITION 3.3\ \rm [G3, Proposition 5.1, p.\ 223]\ \ \it Let $\Psi $ be in the disc 
algebra, let $a\in b\D$ and 
assume that the function  $z\mapsto  \Psi (z)/(z-a)\ (z\in b\D \setminus\{ a\} )$ 
extends continuously to $b\D$. Then there is a function $\Psi _1$ from the disc algebra such that 
$\Psi _1(z)=
\Psi (z)/(z-a)\ (z\in \overline\D \setminus\{ a\})$. 
\vskip 2mm\rm

\noindent\bf Proof of Theorem 3.1 continued. \rm 
Writing $P(z)=p_0(z-a_1)(z-a_2)\cdots (z-a_M)$ we have $M\leq N$ and  
$(z-a_1)\cdots (z-a_M)f(z)=H_1(z)\ (z\in b\D)$ where $H_1 $ belongs to the disc algebra. 
Let $\alpha _1,\cdots ,\alpha _J$ 
be those of $a_1,\cdots, a_M$ which are contained in $\D $. By Proposition 3.3 we may write 
$$ 
(z-\alpha _1)\cdots (z-\alpha _J)f(z) = H(z)\ \ (z\in b\D)
$$ where $H$ belongs to 
the disc algebra and $J\leq N$. This completes the proof of Theorem 3.1.
\vskip 2mm
\noindent \bf REMARK 3.4\ \rm The preceding proof does not work without a smoothness assumption as it is known that 
there are functions $h$ in the disc algebra such that $h(b\D) =h(\overline \D)$ [G1].
\vskip 4mm
\bf 4.\ The general case  \rm
\vskip 2mm
\noindent \bf LEMMA 4.1\ \it Let $N\in\N$ and let $f\colon\ b\D\rightarrow \C$ be a continuous function 
such that the Fourier series 
$$
e^{i(N+1)\theta}f(e^{i\theta})\sim \sum_{n=-\infty }^\infty A_ne^{in\theta}
$$
of the function
$e^{i\theta}\mapsto e^{i(N+1)\theta}f(e^{i\theta})$ is such that $A_0\not=0$ and
$$
A_1=A_2=\cdots =A_N = 0, \ \ A_{-1}=A_{-2}=\cdots=A_{-N} = 0.
$$
There is a polynomial $Q$ such that $f+ Q\not= 0$ on $b\D$  and 
$$
\cW (f+ Q)\leq -N-1 .
$$
\bf Proof.\ \rm With no loss of generality we may assume that $A_0=1$. Since $z^{N+1}f$ is continuous
Fejers theorem 
implies that $z^{N+1}f$ is the uniform limit of the Cezaro means of its Fourier series [Ho]. So, if 
$$
S_k(e^{i\theta})= \sum_{j= -k}^k A_j e^{ij\theta} \ \ k=0, 1, 2,\cdots
$$
are the partial sums of the Fourier series then $e^{i(N+1)\theta}f(e^{i\theta})$ is
the uniform limit, as $m\rightarrow\infty$, of
$$
C_m(e^{i\theta}) = {1\over{m+1}}\bigl[ S_0(e^{i\theta})+S_1(e^{i\theta})+\cdots +S_m(e^{i\theta})\bigr ]
$$
However, each partial sum $S_n$ and therefore each Cezaro mean $C_m$ has the same coefficients 
vanishing property 
as the one which we have assumed for the Fourier series of $z^{N+1}f$:
$$
\hat {(C_m)}(0) = 1,\ \ \hat {(C_m)}(j) = 0\ \ (-N\leq j\leq N,\ j\not= 0)
$$
so that
$$
C_m(z) = 1+z^{N+1}R_m(z) +\overline{z^{N+1}T_m(z)}\ \ (z\in b\D)
$$
where $R_m, T_m$ are polynomials.

Choose $m$ so large that
$$
\vert C_m(z)-z^{N+1}f(z)\vert \leq  {1\over 2}\ \ (z\in b\D) .
\eqno (4.1)
$$
We have
$$
C_m(z)-z^{N+1}R_m(z)-z^{N+1}T_m(z) \in 1 + i\R\ \ (z\in b\D)
$$
which, by (4.1) implies that
$$z^{N+1}f(z)-z^{N+1}R_m(z)-z^{N+1}T_m(z) \in [1/2, 3/2] + i\R\ \ (z\in b\D).
$$ 
It follows that $z^{N+1}f-z^{N+1}R_m-z^{N+1}T_m\not= 0$ on $b\D$ 
and $\cW (z^{N+1}(f-R_m-T_m))= 0$.
Thus 
$$
\cW (f-R_m-T_m)= -N-1
$$ so putting $Q=-R_m-T_m$ completes the proof.
\vskip 2mm
\noindent \bf REMARK 4.2\ \rm Note that the assumption in Lemma 4.1 is equivalent to saying that 
$\hat f(-N-1)\not= 0$ and $\hat f(-1)=\hat f(-2)= \cdots = \hat f(-N) = 0$ 
and $\hat f(-N-2)=\hat f(-N-3)=\cdots =\hat f(-2N-1) = 0$.
\vskip 2mm
We now turn to the proof of Theorem 1.1. We have 
already proved the only if part in Section 1. To prove the if part 
suppose that $f\colon\ b\D\rightarrow\C$ is a continuous function which satisfies (1.1) 
for all polynomials $P,\ Q$ such that $Pf+Q\not=0$ on 
$b\D$. If $N=0$ then 
we already know that $f$ extends holomorphically through $\D$ so assume that 
$N\geq 1$.
\vskip 2mm
\noindent\bf LEMMA 4.3\ \it 
Let $ F\colon \ b\D\rightarrow \C$ be a continuos function. Assume that for some $N\in \N$ we have 
$$
\cW (PF+Q)\geq -N
\eqno (4.2)
$$
whenever $P,\ Q$ are polynomials such 
that $Pf + Q\not= 0$ on $b\D$. Then
$$
F(z) = G(z) + \overline{H(z)} \ \ (z\in b\D )
$$
where $G$ belongs to the disc algebra and $H$ is a rational function 
holomorphic in a neighbourhood of $\overline\D$. \rm
\vskip 2mm
Assume for a moment that Lemma 4.3 holds. 
Since our rational function $H$ is smooth on $b\D$ 
the if part of Theorem 1.1 is now an immediate consequence of Lemma 4.3 and Theorem 3.1.

It remains to prove Lemma 4.3. Given an infinite row $A=(a_1, a_2, \cdots )$ 
and $J\in \N$ we denote by $A(J)$ the row containing the first $J$ entries of $A$,
that is, $A(J)=(a_1, a_2, \cdots ,a_J)$. 

Assume that $F\in C(b\D )$ satisfies (4.2) whenever $P, Q$ are polynomials such 
that $PF+Q\not=0$ on $b\D$. Lemma 4.1 implies that if $P$ is a polynomial such that 
$$
\eqalign{
&\hat{(PF)}(-N-2)=\hat{(PF)}(-N-3)=\cdots = \hat{(PF)}(-2N-1) = 0 \cr
&\hat{(PF)}(-1) = \hat{(PF)}(-2) =\cdots =\hat{(PF)}(-N) = 0 \cr }
$$
then $\hat{(PF)}(-N-1) =0$. If $P(z)=D_0+D_1z+\cdots + D_M z^M$ then
$$
\eqalign{
\hat{(PF)}(-1) &= D_0\hat F(-1)+D_1\hat F(-2)+\cdots +D_M\hat F (-M-1) \cr
&\cdots \cr
\hat{(PF)}(-N+1) &= D_0\hat F(-N+1)+D_1\hat F(-N)+\cdots +D_M\hat F (-M-N+1) \cr
\hat{(PF)}(-N) &= D_0\hat F(-N)+D_1\hat F(-N-1)+\cdots +D_M\hat F (-M-N) \cr
\hat{(PF)}(-N-1) &= D_0\hat F(-N-1)+D_1\hat F(-N-2)+\cdots +D_M\hat F (-M-N-1) \cr
\hat{(PF)}(-N-2) &= D_0\hat F(-N-2)+D_1\hat F(-N-3)+\cdots +D_M\hat F (-M-N-2) \cr
&\cdots \cr
\hat{(PF)}(-2N-1) &= D_0\hat F(-2N-1)+D_1\hat F(-2N-2)+\cdots +D_M\hat F (-2N-M-1). \cr}
$$
Consider the infinite rows 
$$
\eqalign{
&X_{-1}= (\hat F(-1), \hat F(-2), \hat F(-3),\cdots )\cr
&X_{-2}= (\hat F(-2), \hat F(-3), \hat F(-4),\cdots )\cr
&\cdots \cr
&X_{-N}= (\hat F(-N), \hat F(-N-1), \hat F(-N-2),\cdots )\cr
&X_{-N-1}= (\hat F(-N-1), \hat F(-N-2), \hat F(-N-3),\cdots )\cr
&X_{-N-2}= (\hat F(-N-2), \hat F(-N-3), \hat F(-N-4),\cdots )\cr
&\cdots \cr 
&X_{-2N-1}= (\hat F(-2N-1), \hat F(-2N-2), \hat F(-2N-3),\cdots )\cr}
$$
The preceding discussion shows that for every $M\in \N$ the following holds: 
if a row 

\noindent 
$(\overline {D_0}, 
\overline{D_1}\cdots \overline{D_M})
$ is orthogonal to the rows 
$$
\eqalign{
X_{-1}(M+1), X_{-2}(M+1),\cdots ,X_{-N+1}(M+1), X_{-N}(M+1) \cr
X_{-N-2}(M+1), X_{-N-3}(M+1),\cdots ,X_{-2N-1}(M+1) \cr}
\eqno (4.3)
$$
then it is orthogonal to $X_{-N-1}(M+1)$. 
This implies that for every $M\in \N$ the 
row $X_{-N-1}(M+1)$ is a linear 
combination of $2M$ rows (4.3). It follows that there are 
numbers $\lambda _j ,\ -2N-1\leq j\leq -1, \ j\not= -N-1$, such that
$$
X_{-N-1} = \sum_{-2N-1\leq j\leq -1, \ j\not= -N-1}\lambda_j X_j .
\eqno (4.4)
$$

Consider the function 
$$
\Psi (z) = \hat F(-N-1)z +\hat F(-N-2)z^2 +\cdots 
$$ 
The function $\Psi $ is holomorphic on $\D$ and since
$$
\sum_{n=-\infty}^{-N-1}|\hat F(n)|^2 \leq \sum_{n=-\infty}^{\infty}|\hat F(n)|^2 <\infty
$$
it follows that $\Psi $ belongs to the space $H^2$\ [R]. 
We use (4.4) to show that 
$\Psi $ is a rational function. Note that (4.4) implies that
$$
\eqalign{
\Psi (z)  = &\hat F(-N-1)z+\hat F(-N-2)z^2+\cdots  \cr 
= \lambda _{-1}\bigl( &\hat F(-1)z+\hat F(-2)z^2+\cdots \bigr) +\cr
+ \lambda _{-2}\bigl( &\hat F(-2)z+\hat F(-3)z^2+\cdots \bigr) +\cr
&\cdots \cr
+ \lambda _{-N}\bigl( &\hat F(-N)z+\hat F(-N-1)z^2+\cdots \bigr) +\cr
+ \lambda _{-N-2}\bigl( &\hat F(-N-2)z+\hat F(-N-3)z^2+\cdots \bigr) +\cr
&\cdots \cr
+ \lambda _{-2N-1}\bigl( &\hat F(-2N-1)z+\hat F(-2N-2)z^2+\cdots \bigr)\cr}
$$
It follows that 
$$
 \eqalign{\Psi (z)
& = \lambda_{-1}\bigl[P_{-1}(z)+z^N\Psi(z)\bigr] +\cr
& + \lambda_{-2}\bigl[P_{-2}(z)+z^{N-1}\Psi(z)\bigr] +\cr
&  \cdots \cr
& + \lambda_{-N}\bigl[P_{-N}(z)+z\Psi(z)\bigr] +\cr
& + \lambda_{-N-2}z^{-1}\bigl[\Psi (z)-P_{-N-2}(z)\bigr] +\cr
& + \lambda_{-N-3}z^{-2}\bigl[\Psi (z)-P_{-N-3}(z)\bigr] +\cr
&   \cdots \cr
& + \lambda_{-2N-1}z^{-N}\bigl[\Psi (z)-P_{-2N-1}(z)\bigr]\cr}
$$
where 
$$
\eqalign{ 
& P_{-1}(z)=\hat F (-1)z+\hat F(-2)z^2+\cdots + \hat F(-N)z^N \cr
& P_{-2}(z)=\hat F (-2)z+\hat F(-3)z^2+\cdots + \hat F(-N)z^{N-1} \cr
& \cdots\cr
& P_{-N}(z)=\hat F(-N)z\cr
& P_{-N-2}(z) = \hat F (-N-1)z\cr
& P_{-N-3}(z) = \hat F(-N-1)z+\hat F(-N-2)z^2\cr
& \cdots\cr
& P_{-2N-1}=\hat F(-N-1)z+\hat F(-N-2)z^2+\cdots + \hat F (-2N)z^N .\cr}
$$
So
$$\left.\eqalign{
&\Psi(z)\bigl[1-\lambda_{-1}z^N-\lambda_{-2}z^{N-1}-\cdots -
\lambda_{-N}z-\lambda _{-N-2}z^{-1} -\cdots \lambda _{-2N-1}z^{-N}\bigr] = \cr
& =\lambda_{-1}P_{-1}(z)+\lambda _{-2}P_{-2}(z)+\cdots +\lambda_{-N}P_{-N}(z)
-\lambda_{-N-2}P_{-N-2}(z)z^{-1}-\cr
& \ \ \ -  \lambda _{-N-3}P_{-N-3}(z)z^{-2}-\cdots -\lambda _{-2N-1}P_{-2N-1}(z)z^{-N}. \cr}\right\}
\eqno (4.5)
$$
Notice that (4.5) implies that there are polynomials $R, S$ with no common factors such that
$\Psi (z) = R(z)/S(z) $ where S has no zero on $\D$ since $\Psi $ is holomorphic on $\D$.  
 If $\beta _1, \cdots \beta _j$ are those 
poles of $\Psi$ that are contained in $b\D$ then $\Psi^\ast $, the radial limit function 
of $\Psi $, satisfies $\Psi^\ast (z) = \Psi (z)\ (z\in b\D\setminus \{\beta_1,\cdots ,\beta _j\})$. 
However, since $\Psi$ belongs to $H^2$ it follows that
$\Psi|(b\D\setminus \{\beta_1,\cdots ,\beta _j\})$ 
belongs to $L^2(b\D)$ [R] which is impossible if there is a pole on $b\D$ since if 
$e^{i\tau}$ is such a pole then as $\theta\rightarrow \tau$ 
the function $\theta\mapsto  |\Psi (e^{i\theta})|^2$ grows at 
least as fast as a multiple of 
$1/|\theta-\tau|^2$ which is not integrable. Thus, $\Psi $ has no poles on $b\D$ and consequently 
$\Psi $ is a rational function holomorphic in a neighbourhood of $\overline \D$ and so 
$\Phi (z) = \hat F(-1)z+\hat F(-2)z^2 + \cdots = \hat F (-1)z+\cdots+\hat F(-N)z^N+z^N\Psi (z)$ 
is also a rational function holomorphic in a neighbourhood of $\overline\D $. Thus, 
$ H(z) =\overline{\Phi (\overline z)}$ 
is again a rational function holomorphic on a neighbourhood of $\overline \D$. Note that 
$$
\overline {H(z)} = \hat F(-1){1\over z} + \hat F(-2){1\over{z^2}}+\cdots \ \ \ (z\in b\D).
$$
Since $H$ is continuous on $b\D$ it follows that $G=F-\overline H$ is continuous on $b\D$ 
with vanishing Fourier coefficients of negative indices and so $F=G+\overline H $ on $b\D$ where
$G$ is in the disc algebra and $H$ is a rational function holomorphic in a neighbourhood of $\overline\D$.  
This completes the proof of Lemma 4.3. The proof of Theorem 1.1 is complete.
\vskip 4mm
\bf 5.\ Moment conditions and meromorphic extendibility \rm 
\vskip 2mm
Let $f\colon\ b\D\rightarrow \C$ be a continuous function which
extends meromorphically through $\D$ and is such that the meromorphic
extension has at most 
$N$ poles in $\D$, counting mul\-ti\-pli\-ci\-ty. Then there is a nonzero polynomial $P$ of degree
not exceeding $N$ such that 
$Pf$ extends holomorphically through $\D$. Conversely, if $P$ is a nonzero 
polynomial of degree not exceeding $N$ such that $Pf$ extends holomorphically through $\D$ 
then, after using Proposition 3.3 to factor out the zeros of $P$ on $b\D$, we may assume that there 
are a function $H$ in the disc algebra and a polynomial $Q$ of degree not exceeding $N$ 
with all zeros contained in $\D$, such that $f=H/Q$ on $b\D$ which means that $f$ extends 
meromorphically through $\D$ and the meromorphic extension has at most $N$ poles in $\D$. Thus, 
$f$ extends meromorhically through $\D$ with at most $N$ poles, counting multiplicity, 
if and only if there is a nonzero polynomial $P$ of degree not exceeding $N$, such that
$$
\hat {(Pf)}(-n)=0
\eqno (5.1)
$$
for all $n\in\N$.  If $P(z) = D_0+D_1z+D_2z^2+\cdots +D_Nz^N$ then (5.1) means that
$$
D_0\hat f(-n)+D_1\hat f(-n-1)+\cdots +D_N\hat f(-n-N) = 0,
\eqno (5.2)
$$
so  $f$ 
extends meromorphically through $\D$ if and only if there are complex numbers 
$D_0,\cdots ,$ $D_N$, not all zero, such that (5.2) holds for all $n\in\N$. If this happens then the
meromorphic extension of $f$ has at most $N$ poles in $\D$, counting multiplicity. Using the reasoning
applied in Section 3 we can strenghten this to
\vskip 2mm
\noindent \bf PROPOSITION 5.1\ \it Let $f\colon\ b\D\rightarrow \C$ be a continuous function and 
let $N\in\N$. Let $D_0, D_1,\cdots D_N$ be 
a nontrivial solution of the system
$$
\eqalign {
D_0\hat f(-1)+D_1\hat f(-2)+\cdots +D_N\hat f(-1-N) = 0 \cr
D_0\hat f(-2)+D_1\hat f(-3)+\cdots +D_N\hat f(-2-N) = 0 \cr
\cdots \cr
D_0\hat f(-N )+D_1\hat f(-N-1 )+\cdots +D_N\hat f(-2N ) = 0 \cr}
\eqno (5.3)
$$
 The function $f$ has a meromorphic extension through the 
unit disc with at most $N$ poles if and only if these 
numbers $D_0,D_1,\cdots, D_N$ satisfy 
\rm (5.2)\ \it  for all $n\in\N,\ n\geq N+1$. 
\vskip 2mm
\noindent\bf REMARK 5.2\ \rm Note that the system (5.3) is a 
homogeneous system of $N$ linear equations with $N+1$ 
unknowns and so it always has a nontrivial solution.
\vskip 2mm
\noindent\bf Proof of Proposition 5.1.\ \rm Observe first that if 
$a\in\D$ and $k\in\N$ then we have
$$
{1\over{ (z-a)^k}} = \overline{ {z^k}\over {(1-\overline a z)^k}}\ \ \ \ (z\in b\D)
$$
where the function $z\mapsto z^k/(1-\overline az)^k$ is holomorphic in a neighbourhood 
of $\overline\D$. Note also that if $\Phi $ is in the disc algebra, $m\in\N$ and $a\in\D$ then
$$
{ {\Phi (z)} \over {(z-a)^m} } =
{ {\Phi (a)} \over {(z-a)^m} } +
{ {\Phi ^\prime (a)} \over {(z-a)^{m-1}} }+\cdots + { {\Phi^{(m-1)}(a)}\over{(m-1)!(z-a)}} + H_1(z)
$$
where $H_1$ is in the disc algebra. Using decomposition into partial fractions we now 
see that whenever 
$g$ is of the form 
$$
g(z) = { {\Phi (z)} \over {(z-a_1)^{k_1}(z-a_2)^{k_2}\cdots (z-a_J)^{k_J}}    }
$$ 
with $\Phi $ in the disc algebra and $a_j\in\D\ (1\leq j\leq J)$ then
$$
g(z) = F(z) +\overline {G(z)} \ \ \ (z\in b\D)
$$ 
where $F$ is in the disc algebra and $G$ is a rational function 
holomorphic in a neighbourhood of $\overline\D$. 

If $D_0, D_1,\cdots D_N$, not all of them being zero, satisfy (5.3) and (5.2) for all $n\geq N+1$, 
then they satisfy (5.2) for all $n\in N$ so by the preceding discussion $f$ extends meromorphically 
through $\D $ and the meromorphic extension has at most $N$ poles, counting multiplicity. 

To prove the converse, assume that there are numbers $a_1, a_2,\cdots ,a_J$ in $\D$, positive integers 
$k_1, k_2,\cdots , k_J$ such that $k_1+k_2+\cdots +k_J \leq N$, and a function H from the disc algebra 
such that 
$$
f(z) = {{H(z)}\over{(z-a_1)^{k_1}\cdots (z-a_J)^{k_J}}}\ \ (z\in b\D).
$$
By the argument principle it follows that
$$
\left.\eqalign{
\cW (Pf+Q)\geq -N\ \hbox{\ whenever\ } P, Q, \hbox{\ are}\cr
\ \hbox{polynomials such that \ } Pf+Q\not= 0\hbox{\ on\ } b\D .}\right\} 
\eqno (5.4)
$$ 
Let $P(z) = D_0+D_1z+\cdots + D_Nz^N$ be a nonzero polynomial such 
that $\hat{(Pf)}(j) = 0\ \ (-N\leq j\leq -1)$, that is, let $D_0, D_1,\cdots D_N$, not all being zero, 
satisfy (5.3). By the preceding discussion 
$$
(Pf)(z) = F(z) +\overline{ z^{N+1}G(z)} \ (z\in b\D)
$$
where F is in the disc algebra and $G$ is a rational function holomorphic 
in a neighbourhood of $\overline\D $. In particular, 
$G$ is smooth on $b\D$. Assume for a moment that $G\not\equiv 0$. Then, as 
in the proof of Theorem 3.1, we find an $\alpha\in\C$ such that
$Pf-F-\alpha\not=0$ on $b\D$ and that 
$\cW (Pf-F-\alpha )\leq -N-1$. A sufficiently good polynomial 
approximation $Q$ of $-F-\alpha $ then satisfies $\cW (Pf+Q)
\leq -N-1$ which contradicts (5.4). It follows that $G\equiv 0$ 
so $Pf = F$ where $F$ is in the disc algebra and consequently 
$\hat {(Pf)}(j)=0\ (j\leq -N-1)$, that is, (5.2) is satisfied 
for all $n,\ n\geq N+1$. This completes the proof.
\vskip 4mm
\bf 6.\ Remarks\rm
\vskip 2mm
\rm Theorem 1.1 is a one-variable theorem about meromorphic extensions of continuous 
functions on the unit circle. 
It can be described also in more geometric terms as a theorem in $\C\times\overline\C$ as follows.
Let $f\colon\ b\D\rightarrow \C$ be 
a continuous function. Then its graph $\Gamma_f = \{ (z,f(z))\colon\ z\in b\D\}$ is a simple closed 
curve. Suppose that $P,\ Q$ are polynomials such that $Pf+Q\not=0 $ on $b\D$, that is, such that the variety 
$$
V = \{ (z,w)\in\C^2\colon\ P(z)w + Q(z) = 0\} 
\eqno (1.2)
$$
misses $\Gamma_f$. In the special case when $f$ is smooth  then the graf $\Gamma_f$ is 
a smooth curve and the linking number $\hbox{link}(\Gamma_f, V)$ is well defined [AW] 
and is equal to $\cW(Pf+Q)$\ [AW, Lemma 1.2, p.130]. If $f$ is merely continuous then for all 
smooth curves  $\Gamma$ homotopic to 
$\Gamma_f$ in $\C^2\setminus V$ the linking number $\hbox{link}(\Gamma, V)$ is the same which 
implies that for a continuous function $f$ such that $\Gamma _f$ misses $V$ we may 
define $\hbox{link}(\Gamma_f, V)$ simply as 
$\hbox{link}(\Gamma_g, V)$ where $\Gamma_g =\{(z,g(z))\colon\ z\in b\D\}$ 
is the graph of a sufficiently good 
smooth approximation $g\colon\ b\D\rightarrow $ of $f$, so we have $\hbox{link}(\Gamma_f,V)= \cW (Pg+Q) 
= \cW (Pf+Q)$. 

If $f$ has a meromorphic extension $\tilde f$ through $\D $ then the graph $\{(z,\tilde f(z))\colon\ 
z\in\D \}$ of $\tilde f$ is 
a complex submanifold of $\D\times \overline\C$ 
attached to $b\D\times \overline \C$ along $\Gamma_f$. 
So Theorem 1.1 says that the curve $\Gamma _f$ bounds a 
submanifold of $\D\times\overline\C$ 
(that is a graph over $\D$) if and only if the linking numbers $\hbox{link}(\Gamma_f,V)$ for
algebraic varieties $V$ of the form (1.2) which miss $\Gamma_f$, are bounded from below. 
Since $f$ is only assumed to be continuous, our curve $\Gamma_f$ is not smooth,
although, being
a graph over $b\D$ it is quite special. For general curves $\Gamma$ there are recent 
results in a similar spirit, 
 under the assumption of real analyticity and with the assumption 
on linking numbers made for all algebraic varieties which miss $\Gamma $ [HL, Th.6.6].
\vskip 2mm
In Theorem 3.1 it is enough to assume (1.1) 
only for polynomials $P$ of degree not exceeding $N$. We are not able 
to see that the same holds for Theorem 1.1. In our proof of Theorem 1.1 we need polynomials of 
arbitrarily high degree to prove that (1.1) implies that 
$$
f(z)=G(z)+\overline{H(z)}\ \ (z\in b\D)
\eqno (6.1)
$$
where $H$ is a rational function holomorphic in a neighbourhood of $\overline\D$. We can then use 
only the smoothness of $H$ on $b\D$ to be able to apply Theorem 1.1 to show that the 
meromorphic extension of $f$ has 
at most $N$ poles in $\D$. On the other hand, once we know that $f$ is of the form (6.1) where 
$H$ is a rational 
function holomorphic in a neighbourhood of $\overline\D$ then we can, alternatively, show directly 
that the number 
of poles in $\D$ does not exceed $N$ by using
\vskip 2mm
\noindent\bf LEMMA 6.1\ \it Let $\Psi $ be in the disc algebra and let $a_1, a_2\cdots a_m \in\D$ be such that 
$\Psi (a_j)\not= 0\ (1\leq j\leq m)$. There is a polynomial $Q$ such that $\Psi + (z-a_1)\cdots (z-a_m)Q$ 
has no zero on $\overline\D$. \rm 
\vskip 2mm
\noindent Indeed, assuming Lemma 6.1 for a moment, one observes that if $f$ is of 
the form (6.1)where $H$ is a rational function holomorphic 
in a neighbourhood of $\overline\D$ then $f$ must be of the form
$$
f(z) = {{\Psi (z)}\over{(z-a_1)\cdots (z-a_m)}}\ \ (z\in b\D)
$$
where $\Psi$ is in the disc algebra, $a_j\in \D\ (1\leq j\leq m)$ and $\Psi (a_j)\not=0\ (1\leq j\leq m)$. 
Then one uses Lemma 6.1 and (1.1) to show that 
$$
\cW (f+Q)=\cW \Bigl ({{ \Psi +(z-a_1)\cdots (z-a_m)Q}\over{(z-a_1)\cdots (z-a_m)}}\Bigr) = -m\geq -N
$$
so that $m\leq N$ what we wanted to show.
\vskip 2mm
\noindent \bf Proof of Lemma 6.1. \rm 
Let $c_0, c_1,\cdots c_m\in\C,\ c_0\not= 0$. Computing higher order derivatives 
of $z\mapsto e^{\Phi (z)}$ 
it is easy to see that there 
are numbers $d_0, d_1, \cdots, d_m\in\C$ such that if $\Phi $ is holomorphic 
in a neighbourhood of $0$ with Taylor expansion  
$$
\Phi (z) = d_0+d_1z+\cdots d_mz^m+\cdots 
$$
then $e^\Phi$ has Taylor expansion 
$$
e^{\Phi (z) }= c_0 +c_1z+c_2z^2+\cdots c_mz^m +\cdots .
$$

Note that it is enough to construct a $Q$ in the disc algebra as 
then a sufficiently good polynomial approximation of $Q$ will have all the required properties. 

Write $(z-a_1)\cdots (z-a_m)=(z-\alpha_1)^{p_1}\cdots 
(z-\alpha _k)^{p_k}$ where $p_1+\cdots +p_k=m$ and where $\alpha_i\not=\alpha_j\ (i\not= j)$. 
Our $Q$ will have to 
satisfy $\Psi +(z-\alpha_1)^{p_1}\cdots (z-\alpha _k)^{p_k} = e^\Phi$ with $\Phi $ 
from the disc algebra which means that for each $j,\  1\leq j\leq k,$ the function 
$e^\Phi - H$ must have zero of order at least $p_j$ at the point $\alpha _j$. This 
means that for each $j, \ 1\leq j\leq k$, the Taylor expansion of 
$e^\Phi$ in a neighbourhood of $\alpha _j$ has the form
$$
e^{\Phi (z)}=\Psi (\alpha_j)+\cdots + {{\Psi^{(p_j)}(\alpha_j)}
\over{p_j!}}(z-\alpha_j)^{p_j}+\cdots 
\eqno (6.2)
$$
By the preceding discussion there are numbers  $d_{j,\ell},\ 0\leq \ell\leq p_j,\ 1\leq j\leq k$, such that 
if for each $j, \ 1\leq j\leq k$ the function $\Phi $
satisfies
$$
\Phi (z) = d_{j,0}+d_{j,1}(z-\alpha _j)+\cdots +d_{j, p_j}(z-\alpha_j)^{p_j}+\cdots 
$$
then (6.2) holds for each $j,\ 1\leq j\leq k$. It is an easy application of the Weierstrass 
factorization theorem to 
construct an entire function $\Phi $ with this property [R, Th.\ 15.13, p.\ 304]. The function 
$$
Q={{e^\Phi -\Psi}\over{(z-\alpha_1)^{p_1}\cdots (z-\alpha _k)^{p_k}}}
$$
will have the required properties. This completes the proof. 
\vskip 2mm
The following question is open:
\vskip 2mm
\noindent\bf QUESTION 6.2\ \rm Let $f\colon\ b\D\rightarrow \C$ be a continuous function.  
Suppose that for some $N\in\N$ we have 
$\cW (f+Q)\geq -N $ for all polynomials $Q$ such that $f+Q\not= 0$ on $b\D$. 
Must $f$ extend meromorphically through $\D$?
\vskip 2mm \noindent In other words, we are asking whether in Theorem 1.1 it is
enough to assume that $P\equiv 1$ or, equivalently, whether the precise analogue of Theorem 1.0 holds 
for meromorphic extendibility. We do not know the answer even in the case when $f$ is smooth.
\vskip 2mm
We conclude with a remark about holomorphic extendibility. It is an obvious question whether Theorem 1.0 
holds for a smaller class of polynomials $Q$. That linear polynomials do not suffice was shown in [W], 
that polynomials of uniformly bounded degree do not suffice was shown in [G2]. One may ask, for instance, 
whether the polynomials $Q$ satisfying $Q(0)=0$ suffice. The answer is no as shown by the example 
$f(z)=z/(z-1/2)\ (z\in b\D)$. Indeed, writing $Q(z)=zQ_1(z)$ where $Q_1$ is a polynomial 
the argument principle implies that  
$$
\cW (f+Q) = \cW\Bigl( {{z}\over{z-1/2}} +Q \Bigr) 
=\cW \Bigl(  {{z\bigl[1+(z-1/2)Q_1\bigr]}\over {z-1/2}} \Bigr) \geq 0
$$
for all polynomials $Q$ such that $Q(0)=0$  and such that 
$f+Q\not= 0$ on $b\D$, yet $f$ does not extend holomorphically through $\D$. 
However, there is no such example if $f$ has a meromorphic extension through $\D$
which does not vanish at $0$ by the following
\vskip 2mm
\noindent\bf PROPOSITION 6.3\ \it Let $S$ be a polynomial with all its zeros contained in $\D$. 
Suppose that $f\colon\ b\D\rightarrow \C$ is a continuous function such that
$\cW (f+SQ)\geq 0$ for each polynomial $Q$ such that 
$f+SQ\not= 0$ on $b\D$. If $f$ extends meromorphically through $\D$ and 
the meromorphic extension has no 
common zeros with $S$ then $f$ extends holomorphically through $\D$. 
\vskip 2mm
\noindent\bf Proof.\ \rm Suppose that
$$
f(z) = {{H(z)}\over{(z-a_1)\cdots (z-a_N)}}\ \ \ (z\in b\D)
$$
where $a_j\in\D\ (1\leq j\leq N)$ and where $H$ is in the disc algebra,
\ $H(a_j)\not= 0\ (1\leq j\leq N)$, such that $H$ has no common zero with $S$. 
By Lemma 6.1 there is a polynomial $Q$ such that 
$H+(z-a_1)\cdots (z-a_N)SQ $ has no zero on $\overline\D$ It follows that
$$
0\leq \cW (f+SQ)=\cW\Bigl( {{H+(z-a_1)\cdots (z-a_N)SQ}\over{(z-a_1)\cdots (z-a_N)}}\Bigr)= -N
$$
which implies that $N=0$, so $f$ extends holomorphically through $\D$. This completes the proof. 
\vskip 10mm
\noindent\bf ACKNOWLEDGEMENTS \ \rm 
The author is grateful to Sergey Ivashkovich for explaining the simple 
facts about meromorphic extendibility mentioned in the beginning of Section 5.

This work was supported 
in part by the Ministry of Higher Education, Science and Technology of Slovenia 
through the research program Analysis and Geometry, Contract No.\ P1-0291
\vfill
\eject
\vskip 10mm
\centerline{\bf REFERENCES}
\vskip 5mm

\noindent [AW]\ H.\ Alexander and J.\ Wermer: Linking 
numbers and boundaries of varieties. 

\noindent Ann.\ Math.\ 151 (2000) 125-150
\vskip 2mm
\noindent [G1]\ J.\ Globevnik:\ The range of analytic extensions.

\noindent Pacif.\ J.\ Math. 69 (1977) 365-384
\vskip 2mm
\noindent [G2]\ J.\ Globevnik:\ Holomorphic extendibility and the argument principle.

\noindent Complex Analysis and Dynamical Systems II. (Proceedings
  of a conference held in honor of Professor Lawrence Zalcman's sixtieth
  birthday in Nahariya, Israel, June 9-12, 2003), Contemp.\ Math.\ 382 (2005) 171-175
\vskip 2mm
\noindent [G3]\ J.\ Globevnik:\ The argument principle and holomorphic extendibility 
to finite Riemann surfaces.

\noindent  Math.\ Z.\ 253 (2006) 219-225
\vskip 2mm
\noindent [HL]\ F.\ R.\ Harvey and H.\ B.\ Lawson, Jr.: Projective linking and boundaries of 
positive holomorphic chains in projective manifolds, Part I.

\noindent Preprint, http://www.arxiv.org/abs/math.CV/0512379
\vskip 2mm
\noindent [Ho]\ K.\ Hoffman: \it Banach Spaces of Analytic Functions.\rm 

\noindent Prentice Hall, Englewood Cliffs, 1962
\vskip 2mm
\noindent [R]\ W.\ Rudin: \it Real and Complex Analysis.\rm 

\noindent McGraw-Hill, New York, 1987
\vskip 2mm
\noindent [S]\ E.\ L.\ Stout: Boundary values and mapping degree.

\noindent Michig.\ Math.\ J.\ 47 (2000) 353-368
\vskip 2mm
\noindent [W] J.\ Wermer: The argument principle and boundaries of analytic varieties.

\noindent Oper. Theory Adv. Appl., 127, Birkhauser, Basel, 2001, 639-659
\vskip 2mm
\noindent [Z] A.\ Zygmund: \it Trigonometric series. \rm

\noindent Cambridge University Press, Cambridge, New York, 1959

\vskip 10mm
\noindent Institute of Mathematics, Physics and Mechanics

\noindent University of Ljubljana, Ljubljana, Slovenia

\noindent josip.globevnik@fmf.uni-lj.si

\bye